\documentclass[journal,twoside,web]{ieeecolor}
\usepackage{generic}

\makeatletter
\newcommand{\removelatexerror}{\let\@latex@error\@gobble}
\makeatother
\usepackage[T1]{fontenc}
\usepackage[utf8]{inputenc}
\usepackage[english]{babel}
\usepackage{amsmath,amsthm,amssymb}
\usepackage[]{amsfonts}
\usepackage{graphicx}
\usepackage{url}
\usepackage{makecell}
\usepackage[font=footnotesize,labelfont=footnotesize]{caption}
\usepackage{subcaption}

\usepackage{pgfplots}
\pgfplotsset{compat=newest}
\usetikzlibrary{plotmarks}
\usetikzlibrary{arrows.meta}
\usepgfplotslibrary{patchplots}
\usepackage{grffile}
\usepackage{amsmath}

\pgfplotsset{plot coordinates/math parser=false}
  
\usepackage{bbm}

\definecolor{OliveGreen}{rgb}{0, 0.5, 0}

\newtheorem{theorem}{Theorem}[section]
\newtheorem{proposition}[theorem]{Proposition}
\newtheorem{lemma}[theorem]{Lemma}

\theoremstyle{definition}

\newtheorem{remark}[theorem]{Remark}

\newcommand{\EE}{\mathcal{E}}
\newcommand{\GG}{\mathcal{G}}
\newcommand{\HH}{\mathcal{H}}
\newcommand{\PP}{\mathcal{P}}
\newcommand{\VV}{\mathcal{V}}
\newcommand{\WW}{\mathcal{W}}
\newcommand{\FF}{\mathcal{F}}

\newcommand{\NN}{\mathcal{N}}
\newcommand{\II}{\mathcal{I}}

\newcommand{\BB}{\mathcal{B}}

\renewcommand{\SS}{\mathcal{S}}
\newcommand{\TT}{\mathcal{T}}
\newcommand{\XX}{\mathcal{X}}

\newcommand{\CVaR}{\operatorname{CVaR}}

\newcommand{\VI}{\operatorname{VI}}
\newcommand{\SOL}{\operatorname{SOL}}
\newcommand{\argmin}{\operatorname{argmin}}

\newcommand{\real}{\mathbb{R}}
\newcommand{\realnonnegative}{{\mathbb{R}}_{\ge 0}}
\newcommand{\naturalnumbers}{\mathbb{N}}
\newcommand{\setdef}[2]{\{#1 \; | \; #2\}}

\newcommand{\map}[3]{#1:#2 \rightarrow #3}

\newcommand{\Pb}{\mathbb{P}}

\newcommand{\norm}[1]{\ensuremath{\| #1 \|}}
\newcommand{\Bnorm}[1]{\ensuremath{\Big\| #1 \Big\|}}
\newcommand{\until}[1]{[#1]}

\newcommand{\Fhat}{\widehat{F}}
\newcommand{\cl}{\mathrm{cl}}

\newcommand{\HHaff}{\HH_{\mathrm{aff}}}
\newcommand{\HHineq}{\HH_{\mathrm{ineq}}}
\newcommand{\bh}{\bar{h}}

\newcommand{\bhti}{\Pi_{\HHineq} (\bh)}
\newcommand{\hs}{h^*}

\usepackage[normalem]{ulem} 

\newcommand{\oprocendsymbol}{\hbox{$\bullet$}}
\newcommand{\oprocend}{\relax\ifmmode\else\unskip\hfill\fi\oprocendsymbol}
\newcommand{\longthmtitle}[1]{\mbox{}\textup{\textsl{(#1):}}}

\usepackage[prependcaption,colorinlistoftodos]{todonotes}

\title{Stochastic approximation approaches for CVaR-based variational inequalities} 
\author{Jasper Verbree \qquad Ashish Cherukuri \thanks{The authors are with the Engineering and Technology Institute Groningen, University of Groningen. Email: \texttt{\{j.verbree, a.k.cherukuri\}@rug.nl.}}}

\begin{document}

\maketitle
\thispagestyle{empty}
\begin{abstract}
This paper considers variational inequalities (VI) defined by the conditional value-at-risk (CVaR) of uncertain functions and provides three stochastic approximation schemes to solve them. All methods use an empirical estimate of the CVaR at each iteration. The first algorithm constrains the iterates to the feasible set using projection. To overcome the computational burden of projections, the second one handles  inequality and equality constraints defining the feasible set differently. Particularly, projection onto to the affine subspace defined by the equality constraints is achieved by matrix multiplication and inequalities are handled by using penalty functions. Finally, the third algorithm discards projections altogether by introducing multiplier updates. We establish asymptotic convergence of all our schemes to any arbitrary neighborhood of the solution of the VI. A simulation example concerning a network routing game illustrates our theoretical findings.
\end{abstract}

\vspace*{-0.5ex}

\section{Introduction} \label{section introduciton}

\IEEEPARstart{A}{} variety of equilibrium-seeking problems in game theory can be cast as a variational inequality (VI) problem~\cite{FF-JSP:03}. For example, Nash equilibria of a game and Wardrop equilibria of a network routing game both correspond to solutions of a VI under mild conditions. Inspired by real-life, it is natural to perceive the costs or utilities that players wish to optimize in such games as uncertain. Faced with randomness, risk-preferences of players often determines their decisions. Players optimize a risk-measures of random costs in such scenarios and consequently, equilibria can be found by solving a VI involving risk-measures of uncertain costs.  Motivated by this setup, we consider VIs defined by the conditional value-at-risk (CVaR) of random  costs and develop stochastic approximation (SA) schemes to solve them.

\subsubsection*{Literature review}
The most popular way of incorporating uncertainty in VIs is the stochastic variational inequality (SVI) problem, see e.g.~\cite{UVS:13} and references therein. Here, the map associated to the VI is the expectation of a random function. SA methods for solving SVI are well studied~\cite{UVS:13, YC-GL-YO:14}. A key feature of such schemes is the use of an unbiased estimators of the map using any number of sample of the uncertainty. This leads to strong convergence guarantees under a mild set of assumptions. However, the empirical estimator of CVaR, while being consistent, is biased~\cite{RPR-NDS:10}. Therefore, depending on the required level of precision, more samples are required to estimate the CVaR. This biasedness poses challenges in the convergence analysis of SA schemes. For a general discussion on risk-based VIs, including CVaR, and their potential applications, see~\cite{UR:14}. In~\cite{AC:19-cdc}, a sample average approximation method for estimating the solution of CVaR-based VI was discussed. Our work also broadly relates to~\cite{AT-YC-MG-SM:17} and~\cite{CJ-LAP-FM-SM-CS:18} where sample-based methods are used for optimizing the CVaR and other risk measures, respectively. The convergence analysis of our iterative methods consists of approximating the asymptotic behavior of iterates with a trajectory of a continuous-time dynamical system and studying their stability. See~\cite{HJK-DSC:78} and~\cite{VB:08} for a comprehensive account of such analysis.

\subsubsection*{Statement of Contributions}
Our starting point is the definition of the CVaR-based variational inequality (VI), where the map defining the VI consists of components that are the CVaR of random functions. The feasibility set, assumed to be a convex compact set, is defined by a set of inequality and linear equality constraints. We motivate the VI problem using two examples from noncooperative games. Our \emph{first contribution} is the first SA scheme termed as the projected method. This iterative method consists of moving along the empirical estimate of the map defining the VI and projecting each iterate onto the feasibility set. We show that under strict monotonicity, the projected algorithm asymptotically converges to any arbitrary neighborhood of the solution of the VI, where the size of the neighborhood influences the number of samples required to form the empirical estimate in each iteration. Our \emph{second contribution} is the subspace-constrained method that overcomes the computational burden of calculating projections onto the feasibility set by dealing with equality and inequalities differently. In particular, the proximity to satisfying inequality constraints is ensured using penalty functions and iterates are constrained on the subspace generated by linear equality constraints by pre-multiplying the iteration step by an appropriate matrix. We establish that under strict monotonicity, the algorithm converges asymptotically to any neighborhood of the solution of the VI. Our \emph{third contribution} is the multiplier-driven method where projections are discarded altogether by introducing a multiplier for the inequality constraints. The iterates satisfy equality constraints throughout the execution as is the case with the subspace-constrained method by using matrix pre-multiplication. The iterates are shown to converge asymptotically under strict monotonicity to any neighborhood of the solution of the VI. Finally, we demonstrate the behavior of the algorithms using a network routing example. Our work here is an extension of~\cite{JV-AC:20}. Most notably, we have extended the result for the first projection-based algorithm to the case with only a finite number samples of the uncertainty available in each iteration. In addition our second and third algorithm have been modified by introducing a matrix pre-multiplication which constrains iterates to lie in the affine subspace defined by the equality constraints.

\section{Preliminaries} \label{section preliminaries}

\subsubsection{Notation} \label{sec:notation}
Let $\mathbb{R}$ and $\mathbb{N}$ denote the real and natural numbers, respectively. For $N \in \naturalnumbers$, we write ${\until{N} := \{1, 2, \dots, N\}}$. For a given $x \in \mathbb{R}$, we use the notation $[x]^+ := \max(x,0)$. For $x, y \in \real$, the operator $[x]_{y}^+$ equals $x$ if $y > 0$ and it equals $\max\{0,x\}$ if $y \leq 0$. For $x \in \real^n$, we let $x_i$ denote the $i$-th element of $x$, and the $i$-th element of the vector $[x]^+$ is $[x_i]^+$. For vectors $x,y \in \real^n$, $[x]_{y}^+$ denotes the vector whose $i$-th element is $[x_i]_{y_i}^+$, $i \in \until{n}$. The Euclidean norm of $x$ is given by $\norm{x}$. The Euclidean projection of $x$ onto the set $\HH$ is then denoted $\Pi_\HH(x):=\argmin_{y \in \HH} \norm{x-y}$. The $\epsilon$-neigborhood of $x$ is defined as ${\BB_{\epsilon}(x) := \setdef{y \in \real^n}{\norm{y-x} < \epsilon}}$. The closure of a set $\SS \subset \real^n$ is denoted by $\cl(\SS)$. The normal cone to a set $\XX \subseteq \real^n$ at $x \in \XX$ is defined as ${\NN_\XX(x) := \setdef{y \in \real^n}{y^\top(z - x) \leq 0 \enskip \forall z \in \XX}}$. The set ${\TT_\XX(x) := \cl\big(\cup_{y \in \XX} \cup_{\lambda > 0} \lambda(y - x)\big)}$ is referred to as the tangent cone to $\XX$ at $x \in \XX$. 

\subsubsection{Variational inequalities and KKT points}\label{sec:vi-kkt}
For a given map $F: \mathbb{R}^n \rightarrow \mathbb{R}^n$ and a closed set $\HH \subseteq \mathbb{R}^n$, the associated \textit{variational inequality} (VI) problem, $\VI(\HH,F)$, is to find ${h^* \in \HH}$ such that ${(h - h^*)^\top F(h^*) \geq 0}$ holds for all ${h \in \HH}$. The set of all points that solve $\VI(\HH,F)$ is denoted $\SOL(\HH,F)$. The map $F$ is called \textit{monotone} if $\big(F(x) \! - \! F(y)\big)^\top \! (x \! - \! y) \! \ge \! 0$ holds for all ${x,y \in \mathbb{R}^n}$. If the inequality is strict for $x \not = y$, then $F$ is \textit{strictly monotone}. If $\HH$ is nonempty, compact and $F$ is continuous, then $\SOL(\HH,F)$ is nonempty. If $F$ is strictly monotone, then $\VI(\HH,F)$ has at most one solution. Under the \textit{linear independence constraint qualification} (LICQ), we next characterize $\SOL(\HH,F)$ as the set of Karush-Kuhn-Tucker (KKT) points.
\begin{lemma}\longthmtitle{KKT points of $\VI(\HH,F)$}\label{le:KKT}
    Let 
    \begin{equation*}
        \HH := \setdef{h  \in  \real^n }{Ah  =  b,  \enskip q^i(h)  \leq  0, \enskip \forall i \in  [s]},
    \end{equation*}
    where $A \in \real^{l \times n}$, $b \in \real^{l}$ and $l \in \naturalnumbers$, and the functions ${q^i: \mathbb{R}^n \rightarrow \mathbb{R}}$, $i \in \until{s}$, $s \in \naturalnumbers$ are convex and continuously differentiable. For $q(h) \! := \! (q^1(h), \dots, q^s(h))^\top \! \in \! \real^s$, let ${D q(h) \! \in \! \real^{s \times n}}$ be the Jacobian at $h$. For any ${h^* \in \real^n}$, if there exists a multiplier ${(\lambda^*, \mu^*) \in \real^s \times \real^l}$ satisfying
    \begin{equation} \label{KKT}
        \begin{split}
             &F(h^*) + \big( D q(h^*)\big)^\top \lambda^* + A^\top \mu^* = 0, \\
             Ah^* & = b, \quad q(h^*) \leq 0, \quad  \lambda^{*} \geq 0, \quad \lambda^{*\top}q(h^*) = 0,
        \end{split}
    \end{equation}
    then we have $h^* \in \SOL(\HH,F)$. Such a point $(h^*,\lambda^*, \mu^*)$ is referred to as a KKT point of the $\VI(\HH,F)$. Conversely, for $h^* \in \SOL(\HH,F)$, let $\II_{h^*} = \setdef{i \in \until{s}}{q^i(h^*) = 0}$. If the vectors $\{\nabla q^i(h^*)\}_{i \in \II_{h^*}}$ and the row vectors $\{A_j\}_{j \in [l]}$ are linearly independent, or in other words, the LICQ holds at $h^*$, then there exists a $(\lambda^*, \mu^*)$ satisfying \eqref{KKT}.
\end{lemma}
The above result is well known in the context of convex optimization. The extension to the VI setting can be deduced from $\text{\cite[Proposition 3.46]{DS:18}}$, \cite[Theorem 12.1]{UF-WK-GS:10}, and noting that if $h^* \in \SOL(\HH, F)$, then it is also a minimizer of the function $y \mapsto y^\top F(h^*)$ subject to $y \in \HH$.

\subsubsection{Projected dynamical systems}\label{sec:proj}
For a given map ${F: \mathbb{R}^n \times [0,\infty) \rightarrow \mathbb{R}^n}$ and a closed set $\HH \subseteq \mathbb{R}^n$ the associated projected dynamical system is given by
\begin{equation*}
	\dot{h}(t) = \Pi_{\TT_{\HH}(h(t))}\big(F(h,t)\big).
\end{equation*}
Here $\TT_{\HH}(h)$ is the tangent cone of $\HH$ at $h$ ${\text{(see Section~\ref{sec:notation})}}$.
We say that a map $\bh : [0,\infty) \rightarrow \HH$ with $\bh(0) \in \HH$ is a solution of the above system when $\bh(\cdot)$ is absolutely continuous and $\dot{\bh}(t) = \Pi_{\TT_{\HH}(\bh(t))}\Big(F\big(\bh(t), t\big)\Big)$ for almost all $t \in [0, \infty)$. Note that $\bh(t) \in \HH$ for all $t$. We use the terms solution and trajectory interchangeably throughout the paper. 

\subsubsection{Conditional Value-at-Risk}
The \textit{Conditional Value-at-Risk} (CVaR) \textit{at level} $\alpha \in (0,1]$ of a real-valued random variable $Z$, defined on a probability space $(\Omega, \FF, \mathbb{P})$, is
\begin{equation*} 
    \CVaR_\alpha[Z] := \inf_{\eta \in \mathbb{R}}\big\{\eta + \alpha^{-1} \mathbb{E}[Z - \eta]^+\big\},
\end{equation*}
where the expectation is with respect to $\Pb$. The value $\alpha$ is a constant that characterizes risk-averseness. 
Given $N$ i.i.d samples $\{\widehat{Z}_j\}_{j \in [N]}$ of random variable $Z$, one can approximate $\CVaR_\alpha[Z]$ using the following empirical estimate:
\begin{equation} \label{estimator CVaR}
    \widehat{\CVaR}_\alpha^N[Z] = \inf_{\eta \in \mathbb{R}}\big\{ \eta + (N \alpha)^{-1}  \textstyle\sum_{j=1}^N[\widehat{Z}_j - \eta]^+ \big\}.
\end{equation}
This estimator is biased, but consistent \cite[Page 300]{AS-DD-AR:14}. That is, $\widehat{\CVaR}_\alpha^N[Z] \to \CVaR_\alpha[Z]$ almost surely as $N \to \infty$. 

\section{Problem statement and motivating examples} \label{section problem statement}
Consider a set of functions $C_i: \mathbb{R}^n \times \mathbb{R}^m \rightarrow \mathbb{R}$, ${i \in \until{n}}$, $(h,\xi) \mapsto C_i(h,\xi)$, where $\xi$ represents a random variable with distribution $\Pb$. For a fixed $h$, $C_i(h,\xi)$ is therefore a real-valued random variable. Define the map $\map{F_i}{\real^n}{\real}$ as the CVaR of $C_i$ at level $\alpha \in (0,1]$:
\begin{equation} \label{definition elements F}
       F_i(h) := \CVaR_\alpha\big[C_i(h,\xi)\big], \quad \text{ for all } i \in \until{n}.
\end{equation}
For notational convenience, let $C: \mathbb{R}^n \times \mathbb{R}^m \rightarrow \mathbb{R}^n$ and $\map{F}{\real^n}{\real^n}$ be the element-wise concatenation of the maps $\{C_i\}_{i \in \until{n}}$ and $\{F_i\}_{i \in \until{n}}$, respectively. Let $\HH \subseteq \real^n$ be a nonempty closed set of the form
\begin{equation} \label{explicit contrained set 2}
    \HH := \setdef{h  \in  \real^n }{Ah  =  b,  \enskip q^i(h)  \leq  0, \enskip \forall i  \in  [s]},
\end{equation}
where $A \in \real^{l \times n}$, $b \in \real^l$ and $l \in \naturalnumbers$, and the functions ${q^i: \mathbb{R}^n \rightarrow \mathbb{R}}$, ${i \in \until{s}}$, $s \in \naturalnumbers$ are convex and continuously differentiable. The objective of this paper is to provide stochastic approximation (SA) algorithms to solve the variational inequality problem $\VI(\HH,F)$. Our strategy will be to use an empirical estimator, derived from samples of $C(h,\xi)$, of the map $F$ at each iteration of the algorithm. Below we discuss two motivating examples for our setup.

\subsection{CVaR-based routing games} \label{section routing game}
Consider a directed graph $\GG = (\VV,\EE)$, where ${\VV = \until{n}}$ is the set of vertices, and ${\EE \subseteq \VV \times \VV}$ is the set edges. To such a graph we associate a set $\WW \subseteq \VV \times \VV$ of origin-destination (OD) pairs. An OD-pair $w$ is given by an ordered pair $(v^w_o,v^w_d)$, where $v^w_o, v^w_d \in \VV$ are called the origin and the destination of $w$, respectively. The set of all paths in $\GG$ from the origin to the destination of $w$ is denoted $\PP_w$. The set of all paths is given by $\PP = \cup_{w \in \WW} \PP_w$. Each of the participants, or agents, of the routing game is associated to an OD-pair, and can choose which path to take to travel from its origin to its destination. The choices of all agents give rise to a flow vector $h \in \mathbb{R^{|\PP|}}$. We consider a non-atomic routing game and so $h$ is a continuous variable.

For each (OD)-pair $w$, a value $D_w \geq 0$ defines the demand associated to it. The feasible set $\HH \subset \real^{|\PP|}$ is then given by ${\HH = \setdef{h}{\sum_{p \in \PP_w}h_p = D_w, \enskip \forall w \in \WW, \enskip h_p \geq 0 \enskip \forall p \in \PP }}$. To each of the paths $p \in \PP$, we associate a cost function $C_p: \mathbb{R}^{|\PP|} \times \mathbb{R}^{m} \rightarrow \mathbb{R}, (h,\xi) \mapsto C_p(h,\xi)$, which depends on the flow $h$, as well as on the uncertainty $\xi \in \mathbb{R}^m$. Each agent chooses $p \in \PP_w$ that minimizes $\CVaR_\alpha\big[C_p(h,\xi)\big]$. These elements define the CVaR-based routing game~\cite{AC:19-cdc} to which we assign the following notion of equilibrium: the flow $h^* \in \HH$ is said to be a CVaR-based Wardrop equilibrium (CWE) of the CVaR-based routing game if, for all $w \in \WW$ and all $p,p' \in \PP_w$ such that $h^*_p > 0$, we have
\begin{equation*}
    \CVaR_\alpha\big[C_p(h^*,\xi)\big] \leq \CVaR_\alpha\big[C_{p'}(h^*,\xi)].
\end{equation*}
The intuition is that at equilibrium, for each agent, there is no other path than the selected one that has a smaller value of conditional value-at-risk. Under continuity of $C_p$, the set of CWE is equal to the set of solutions of $\VI(\HH,F)$, where $F: \mathbb{R}^{|\PP|} \rightarrow \mathbb{R}^{|\PP|}$ takes the form~\eqref{definition elements F}.

\subsection{CVaR-based Nash equilibrium}
A more general example of our setup would be in finding the Nash equilibrium of a non-cooperative game~\cite[section 1.4.2]{FF-JSP:03}. Let there be $N$ players with individual cost functions $\map{\theta_i}{\real^{nN}}{\real}$, $x \mapsto \theta_i(x)$ and possible strategy sets $\mathcal{X}_i \subseteq \real^n$. Here $x \in \real^{nN}$ denotes the vector containing the strategies of all players, where $x_i \in \real^n$ is the strategy of player $i$. We assume without loss of generality that the action/strategy sets of each player are of the same dimension $n$. An alternative notation is $\theta_i(x) = \theta_i(x_i, x_{-i})$, where $x_{-i}$ is the vector containing the strategies of all players except $i$. Each player $i$ aims to minimize its cost $\theta_i$ by choosing its own strategy optimally. That is, for any fixed $\tilde{x}_{-i}$ they solve
\begin{align*}
	\text{minimize }    \quad &\theta_i(x_i,\tilde{x}_{-i}),  \\
	\text{subject to }  \quad &x_i \in \mathcal{X}_i.
\end{align*}
A Nash equilibrium of such a game is a solution vector $x^*$ such that none of the players can reduce their cost by changing their strategy. Under the assumption that the sets $\mathcal{X}_i$ are convex and closed, and the functions $x_i \mapsto \theta_i(x_i, \tilde{x}_{-i})$ are convex and continuously differentiable for any $\tilde{x}_{-i}$, a joint strategy vector $x^*$ is a Nash equilibrium if and only if it is a solution to $\VI(\mathcal{X}, F)$, where $F(x) := (\nabla_{x_i}\theta_i(x))_{i = 1}^N$ is the concatenation of the gradients of $\theta_i$ functions, and ${\mathcal{X} = \prod_{i = 1}^n\mathcal{X}_i}$ . Consider the functions $\theta_i$ of the form
\begin{equation*}
	\theta_i(x) := \CVaR_\alpha [f_i(x_i, x_{-i})g(\xi) + \bar{f}_i(x_i, x_{-i})],
\end{equation*}
where functions $f_i$, $g$ and $\bar{f}_i$ are real-valued, $\xi$ models the uncertainty, and $f_i(x_i,x_{-i}) \geq 0$ for all $x$. Then,  $\VI(\mathcal{X}, F)$ is a CVaR-based variational inequality as discussed in our paper. Specifically, in this case, since CVaR is positive-homogeneous and shift-invariant~\cite[Chapter 6]{AS-DD-AR:14}, we have
\begin{align*}
	\theta_i(x) = \CVaR_\alpha[g(\xi)] f_i(x_i,x_{-i}) + \bar{f}_i(x_i, x_{-i}).
\end{align*}
As a consequence, we get
\begin{align*}
	\nabla \theta_i(x) = \CVaR_\alpha[g(\xi)] \nabla_{x_i} f_i(x_i,x_{-i}) + \nabla_{x_i} \bar{f}_i(x_i, x_{-i}).
\end{align*}
Under the assumption that $\nabla_{x_i} f_i$ is nonnegative for all ${x \in \mathcal{X}}$, we get
\begin{align*}
	\nabla \theta_i(x) = \CVaR_\alpha[g(\xi) \nabla_{x_i} f_i(x_i,x_{-i}) + \nabla_{x_i} \bar{f}_i(x_i, x_{-i})].
\end{align*}
where $\CVaR$ is understood component-wise. Thus, $F$ can be written as concatenation of $\CVaR$ of various functions and finding the Nash equilibrium of this game is equivalent to solving $\VI(\mathcal{X}, F)$, which fits into our presented framework.

\section{Stochastic approximation algorithms for solving $\VI(\HH, F)$} \label{section results} \label{section stochastic approximation algorithms}
In this section, we introduce the SA algorithms along with their convergence analysis. All introduced schemes approximate $F$ with the estimator given in~\eqref{estimator CVaR}. Given $N$ independently and identically distributed samples $\big\{(\widehat{C_i(h,\xi)} )_j\big\}_{j=1}^N$ of the random variable $C_i(h,\xi)$, let 
\begin{align*}\label{eq:Fhat}
    \Fhat^N_i(h) := \textstyle \inf_{t \in \mathbb{R}}\Big\{ t + (N \alpha)^{-1} \sum_{j=1}^N\big[(\widehat{C_i(h,\xi)})_j - t\big]^+ \Big\}
\end{align*}
stand for the estimator of $F_i(h)$. Analogously, the estimator of $F(h)$ formed using the element-wise concatenation of $\Fhat^N_i(h)$, $i \in \until{n}$, is denoted by $\Fhat^N(h)$. We assume that the $N$ samples of each cost function are a result of the same set of $N$ events, that is, the distribution of $\Fhat^N(h)$ depends on $\Pb^N$. We next present our first algorithm.

\subsection{Projected algorithm}
For a given sequence of step-sizes $\{\gamma^k\}_{k=0}^\infty$, with $\gamma^k > 0$ for all $k$, a sequence $\{N_k\}_{k=0}^\infty \subset \naturalnumbers$, and an initial vector $h^0 \in \HH$, the first algorithm under consideration, which we will refer to as the \textit{projected algorithm}, is given by
\begin{equation}\label{eq:projected-original-form}
    h^{k+1} = \Pi_{\HH}\big( h^k - \gamma^k\Fhat^{N_k}(h^{k})\big),
\end{equation}
where $\Pi_{\HH}$ is the projection operator (see Section~\ref{sec:notation}) and $h^k$ is the $k$-th iterate of $h$ produced by the algorithm. The above algorithm is inspired by the SA schemes for solving a stochastic VI problem, see~\cite{UVS:13} for details on other SA schemes. The key difference from the setup in~\cite{UVS:13} is the fact that there the map $F$ is the expected value of a random variable for which an unbiased estimator $\widehat{F}$ is available. In our case the estimator is biased posing limitations on the sample requirements for convergence of the algorithms. We can write the projected algorithm~\eqref{eq:projected-original-form} equivalently as
\begin{equation} \label{algorithm with error term}
    h^{k+1} = \Pi_{\HH}\Big( h^k - \gamma^k\big(F(h^{k}) + \widehat{\beta}^{N_k} \big) \Big),
\end{equation}
where $\widehat{\beta}^{N_k} := \widehat{F}^{N_k}(h^k) - F(h^k)$ is the error introduced by the estimation. For this and the upcoming algorithms, common assumptions on the sequence $\{\gamma^k\}$ are
\begin{equation} \label{stepsize}
    \begin{split}
       \textstyle \sum_{k = 0}^\infty \gamma^k = \infty, &\qquad  \textstyle \sum_{k = 0}^\infty (\gamma^{k})^2 < \infty.
    \end{split}
\end{equation}
Our first result gives sufficient conditions for convergence of \eqref{algorithm with error term} to any neighborhood of the solution $h^*$ of $\VI(\HH,F)$.
\begin{proposition} \longthmtitle{Convergence of the projected algorithm~\eqref{eq:projected-original-form}} \label{projected theorem}
	Let $F$ as defined in \eqref{definition elements F} be a strictly monotone, continuous function, and let $\HH$ be a compact convex set of the form~\eqref{explicit contrained set 2}. For the algorithm \eqref{algorithm with error term}, assume that step-sizes $\{\gamma^k\}$ satisfies \eqref{stepsize} and the sequence $\{N_k\}$ is such that $\{\widehat{\beta}^{N_k}\}$ is bounded with probability one. Then, for any $\epsilon > 0$ there exists $N_{\epsilon} \in \mathbb{N}$ such that $N_k \geq N_\epsilon$ for all $k$ implies, with probability one,
	\begin{equation*}
		\lim_{k \rightarrow \infty}\norm{h^k - h^*} \le \epsilon.
	\end{equation*}
\end{proposition}
\begin{proof}
	To ease the exposition of this proof, we split the error as $\widehat{\beta}^{N_k} = e^{N_k} + \widehat{\varepsilon}^{N_k}$, where $e^{N_k} = \mathbb{E}[\widehat{\beta}^{N_k}]$. Note that we then have $\mathbb{E}[\widehat{\varepsilon}^{N_k}] = 0$, and  by the boundedness assumption, there exists a constant $B_e > 0$ such that $\norm{e^{N_k}} \leq B_e$ for all $k$. The first step of the proof is to show that the sequence $\{h^k\}$ converges to a trajectory of the following continuous-time projected dynamical system:
	\begin{equation} \label{eq:ODE projected} 
		\dot{\bar{h}}(t) = \Pi_{\TT_\HH(\bar{h}(t))} \Big( -F\big(\bar{h}(t)\big) - e(t) \Big), \quad \bh(0) \in \HH. 
	\end{equation}
	Here $e(\cdot)$ is a uniformly bounded measurable map satisfying $\norm{e(t)} \leq B_e$ for all $t$ (see Section~\ref{sec:proj} for further details on how solutions to projected dynamical systems are defined). For the sake of rigor, we note that the existence of a trajectory of~\eqref{eq:ODE projected} starting from any point in $\HH$ is guaranteed by~\cite[Lemma A.1]{MS-AC-NM:22-auto}. To make precise the convergence of the sequence generated by~\eqref{algorithm with error term} to a trajectory of~\eqref{eq:ODE projected}, we say that $\{h^k\}$ converges to a trajectory $\bar{h}( \cdot )$ of \eqref{eq:ODE projected} if
	\begin{equation} \label{eq:convergence to trajectory}
		\lim_{i \rightarrow \infty} \sup_{j \geq i}\Bnorm{h^j - \bar{h}\Bigl(\sum_{k=i}^{j-1} \gamma^k\Bigr)} = 0.
	\end{equation}
    That is, the discrete-time trajectory formed by the linear interpolation of the iterates $\{h^k\}$ approaches the continuous time trajectory $t \mapsto \bar{h}(t)$. The proof of the existence of a map $\bar{h}(\cdot)$ satisfying~\eqref{eq:convergence to trajectory} is similar to that of~\cite[Theorem 5.3.1]{HJK-DSC:78}, with the only change being the existence of an error term $e(t)$ in dynamics~\eqref{eq:ODE projected} which is absent in the cited reference. The inclusion of the error term is facilitated by reasoning presented in the proof of~\cite[Theorem 5.2.2]{HJK-DSC:78}. We avoid repeating these arguments here in the interest of space. 

    Convergence of the sequence $\{h^k\}$ can now be analyzed by studying the asymptotic stability of \eqref{eq:ODE projected}. To this end, we consider the candidate Lyapunov function
    \begin{equation*} 
	    V\big(\bar{h}\big) = \frac{1}{2}\norm{\bar{h} - h^*}^2,
    \end{equation*}
    where $h^*$ is the unique solution of $\VI(\HH,F)$. We first look at the case $e(\cdot) \equiv 0$.  For notational convenience, define the right-hand side of~\eqref{eq:ODE projected} in such a case by the map ${\map{X_{e\equiv0}}{\real^n}{\real^n}}$. The Lie derivative of $V$ along $X_{e\equiv 0}$ is then given by
    \begin{align} \label{eq:Lie-derivative}
	    \nabla V(\bar{h})^\top X_{e \equiv 0}(\bar{h})  = (\bar{h} - h^*)^\top \Pi_{\TT_\HH(\bar{h})} \big( -F(\bar{h}) \big).
    \end{align}
    We want to show that the right-hand side of the above equation is negative for all $\bar{h} \neq h^*$. We first note that by Moreau’s decomposition theorem~\cite[Theorem 3.2.5]{JBHU-CL:93}, for any $v \in \real^n$ and $\bar{h} \in \HH$, we have ${\Pi_{\TT_\HH(\bar{h})}(v) = v - \Pi_{\NN_\HH(\bar{h})}(v)}$, where $\NN_\HH(\bar{h})$ is the normal cone to $\HH$ at $\bar{h}$. Using the above relation in~\eqref{eq:Lie-derivative} gives
    \begin{align}
	    \nabla V(\bar{h})^\top \! X_{e \equiv 0}(\bar{h}) \! =  &-(\bar{h}  -  h^*)^\top \! F(\bar{h}) \notag   \\
	    &+ (h^*  -  \bar{h})^\top \! \Pi_{\NN_\HH(\bar{h})} \big( \! - \! F(\bar{h}) \big)	\notag \\
	    \leq &-(\bar{h}  -  h^*)^\top \! F(\bar{h}), \label{eq:lyap-ineq-1}
    \end{align}
    where the inequality is due to the definition of the normal cone (see Section \ref{sec:notation}) and the fact that $h^* \in \HH$. 
    Due to strict monotonicity of $F$, we have ${(\bar{h} - h^*)^\top F(\bar{h}) > (\bar{h} - h^*)^\top F(h^*)}$ whenever ${\bar{h} \neq h^*}$. Furthermore, since $h^* \in \SOL(\HH,F)$ we know that ${(\bar{h} - h^*)^\top F(h^*) \geq 0}$ for all $\bar{h} \in \HH$. Combining these two facts implies that the function $W(\bar{h}) := (\bar{h} - h^*)^\top F(\bar{h})$ satisfies $W(\bar{h}) > 0$ whenever $\bar{h} \neq h^*$.
    Using this in the inequality~\eqref{eq:lyap-ineq-1} yields
    \begin{equation} \label{eq:Lyapunov-less-than}
    	\nabla V(\bar{h})^\top X_{e \equiv 0}(\bar{h}) \le - W(\bar{h}) < 0
    \end{equation}
    whenever $\bar{h} \neq h^*$. Now let $\overline{\HH}_{\epsilon} := \setdef{h \in \HH}{\norm{h - h^*} \geq \epsilon}$. Since $\HH$ is compact, $\overline{\HH}_{\epsilon}$ is compact. Since $W$ is continuous, there exists a $\delta > 0$ such that $W(\bar{h}) \ge \delta$ for all $\bar{h} \in \overline{\HH}_{\epsilon}$. Therefore we get, from~\eqref{eq:Lyapunov-less-than},
    \begin{align}\label{eq:Lie-delta}
    	\nabla V(\bar{h})^\top X_{e \equiv 0}(\bar{h}) \le - \delta, \quad \text{for all } \bar{h} \in \overline{\HH}_{\epsilon}.
    \end{align}
    Next, we drop the assumption that $e(\cdot) \equiv 0$ and use the map $\map{X}{\real^n \times [0,\infty)}{\real^n}$ to denote the right-hand side of~\eqref{eq:ODE projected}. Consider any trajectory $t \mapsto \bar{h}(t)$ of~\eqref{eq:ODE projected}. Since the map is absolutely continuous and $V$ is differentiable, we have for almost all $t \ge 0$ and for $\bar{h}(t) \in \overline{\HH}_\epsilon$,
    \begin{align*}
    	\frac{dV}{dt}(t) = \nabla V(\bar{h}(t))^\top X(\bar{h}(t),t) \le - \delta - (\bar{h}(t) - h^*)^\top e(t),
    \end{align*}
    where for obtaining the above inequality we have first used Moreau's decomposition as before to get rid of the projection operator in $X$ and then employed~\eqref{eq:Lie-delta}. Next we bound the error term in the above inequality. Since $\HH$ is compact and $ h^* \in \HH$, there exists $B_h > 0$ such that $\norm{\bar{h} - h^*} \le B_h$ for all $\bar{h} \in \HH$.  In addition $\norm{e(t)} \le B_e$ for all $t$, where $B_e$ is the bound satisfying $\norm{e^{N_k}} \le B_e$. 
    Since the empirical estimate of the CVaR is consistent, we know that $B_e$ can be made arbitrarily small by selecting $N_k$ to be appropriately large for all $k$. That is, there exists $N_\epsilon \in \mathbb{N}$ such that when $N_k > N_\epsilon$ we have $\norm{e^{N_k}} < \frac{\delta}{B_h}$. Consequently, if $N_k > N_\epsilon$ for all $k$, then $\norm{e(t)} < \frac{\delta}{B_h}$ for all $t$. By selecting such a sample size at each iteration and thus bounding the error term, we obtain
    \begin{equation*}
        \begin{split}
            \frac{dV}{dt}(t)&\leq -\delta - (\bar{h}(t) - h^*)^\top e(t)    \\
        	&\leq -\delta + \norm{\bar{h}(t) - h^*} \norm{e(t)} < -\delta + B_h \frac{\delta}{B_h} \leq 0,	
        \end{split}
    \end{equation*}
    which holds for almost all $t$ whenever $\bar{h}(t) \in \overline{\HH}_\epsilon$. That is, the trajectory converges to the set $\setdef{h \in \HH}{\norm{h - h^*} \le \epsilon}$ as $t \to \infty$. This concludes the proof. 
\end{proof}
In the above result, the restriction $N_k \geq N_\epsilon$ does not need to hold for all $k$. The result also holds if there exists a $K \in \naturalnumbers$ such that $N_k \geq N_\epsilon$ for all $k \geq K$. Regarding boundedness of  $\{\widehat{\beta}^{N_k}\}$, it is ensured if for example each $C_i$ is bounded over the set  $\HH \times \Xi$, where $\Xi$ is the support of  $\xi$. 

Despite the convergence property established in Proposition \ref{projected theorem}, the algorithm in \eqref{algorithm with error term} suffers from some disadvantages. Most notably, the algorithm requires computing projections onto the set $\HH$ at each iteration, which can be computationally expensive. To address these issues we propose two algorithms that achieve similar convergence to any neighborhood of the solution of the $\VI(\HH,F)$. The first requires projection onto inequality constraints only and the second does not involve any projection on the primal iterates and instead ensures feasibility using dual variables. As in Proposition~\ref{projected theorem}, we will impose continuity and monotonicity assumptions on $F$ in the upcoming results. We provide the following general result on the continuity and monotonicity properties of $F$.
\begin{lemma}\longthmtitle{Sufficient conditions for monotonicity and continuity of $F$}
    The following hold:
    \begin{itemize}
        \item   If for any $\epsilon > 0$ there exist a $\delta > 0$ such that ${\norm{h - h'} \le         \delta}$ implies $\norm{ C_i(h,\xi) - C_i(h',\xi)} \le \epsilon$ for all $i \in           \until{n}$ and all $\xi$, then $F$ is continuous.
        \item   Assume that there exist functions ${f_i:\real^n \rightarrow \real}$ and ${g_i:            \real^m \rightarrow \real}$ such that ${C_i(h,\xi) \equiv f_i(h) + g_i(\xi)}$, for         all $i \in [n]$. Let ${f(h):= \big(f_1(h), \dots, f_n(h)\big)}$. Then, $F$ is             monotone (resp. strictly) if $f$ is monotone (resp. strictly monotone).
    \end{itemize}
\end{lemma} 
\begin{proof}
    Continuity follows by arguments similar to the proof of \cite[Lemma IV.8]{AC:19-cdc}. For the second part, note that CVaR satisfies ${\CVaR_\alpha\big[C_i(h,\xi)\big] = f_i(h) + \CVaR_\alpha\big[g_i(\xi)\big]}$, for all $h$ and $i \in [n]$ \cite[Page 261]{AS-DD-AR:14}. The proof then follows from the fact that $F(h) - F(h^*) = f(h) - f(h^*)$.
\end{proof}
In the above result, the continuity condition, that may be difficult to check in practice, holds if $\xi$ has a compact support and for any fixed $\xi$, the functions $C_i$ are continuous with respect to $h$. We now introduce our next algorithm.

\subsection{Subspace-constrained algorithm}
In this section, we take a closer look at the form of $\HH$ given in~\eqref{explicit contrained set 2} and design an algorithm that handles inequality and equality constraints independently. To this end, we write $\HHaff  := \setdef{h \in \real^n}{Ah=b}$, and ${\HHineq  := \setdef{h \in \real^n}{q^i(h) \le 0, \enskip \forall i \in [s]}}$ for the sets of points satisfying the equality and inequality constraints, respectively. We then have $\HH = \HHaff \cap \HHineq$. It turns out that, using matrix operation, we can ensure that the iterates of our algorithm always remain in $\HHaff$. The method works as follows. Let $\{a_1,\cdots,a_l\}$ be the row vectors of $A$, and let $\{u_1,\cdots,u_n\}$ be an orthonormal basis for $\mathbb{R}^n$ such that the first $M \in \naturalnumbers$ vectors $\{u_1,\cdots,u_{M}\}$ form a basis for the span of vectors $\{a_1,\cdots,a_l\}$. Then, for the subspace ${\SS = \setdef{g \in \real^n}{Ag = 0}}$,  we have
\begin{equation*}
	\Pi_{\SS}(v) = \textstyle \Big(I - \sum_{i = 1}^{M} u_i u_i^\top \Big)v, \quad \text{ for any } v \in \real^n.
\end{equation*}
This well known fact follows from \cite[Theorem 7.10]{DCL-SRL-JJM:16} and noting that $\Pi_{\SS}(v) = v - \Pi_{\SS^{\perp}}(v)$, where $\SS^{\perp}$ is the set of vectors orthogonal to the subspace $\SS$. Thus, the projection onto $\SS$ is achieved by pre-multiplying with the matrix
\begin{equation} \label{eq:subspace-projection-matrix}
	L := I - \sum_{i = 1}^{M} u_i u_i^T.
\end{equation}
Consequently, for any vector $z$ of the form $z = L v$, ${v \in \real^n}$ we have $A z = 0$. To construct $L$ one can find the orthonormal basis vectors $\{u_i\}_{i \in [l]}$ for the span of $\{a_j\}_{j \in [l]}$ and $\real^n$ by using Gram-Scmhidt orthogonalization process \cite[Section 6.4]{DCL-SRL-JJM:16}. Alternatively, if $A$ has full row rank one can use ${L := I - A^\top(AA^\top)^{-1}A}$, see e.g.,~\cite{CDM:00}. We use this projection operator to define our next method called the \textit{subspace-constrained algorithm}:
\begin{equation}\label{eq:subspace-algo}
	h^{k+1} \! = \! h^k \! - \! \gamma^k L \Big(F(h^k) \! + \! c\big(h^k \! - \! \Pi_{\HH_{\mathrm{ineq}}}(h^k)\big)  \!+ \! \widehat{\beta}^{N_k}\Big),
\end{equation}
where the initial iterate $h^0 \in \HHaff$. In the above, $c > 0$ is a parameter to be specified later in the convergence result, the error sequence $\{\widehat{\beta}^{N_k}\}$ is as defined in~\eqref{algorithm with error term}, and $L \in \real^{n \times n}$ is as defined in \eqref{eq:subspace-projection-matrix}.
	
Due to the presence of $L$ in the above algorithm, the direction in which the iterate moves in each iteration is projected onto the subspace $\SS$. Hence, $h^k \in \HHaff$ for all $k$. We formally establish this in the below result. Furthermore, convergence to a neighbourhood of the set $\HHineq$ is achieved through the term $h^k - \Pi_{\HHineq}(h^k)$. That is, the higher the value of the design parameter $c$, the closer the limit of $\{h^k\}$ is to $\HHineq$. Together, these mechanisms ensure that we keep iterates close to $\HH$ and ultimately drive them to a neighbourhood of $\hs$.
\begin{proposition} \label{pr:subspace-cons} \longthmtitle{Convergence of the subspace-constrained algorithm~\eqref{eq:subspace-algo}}
	Let $F$ as defined in \eqref{definition elements F} be a strictly monotone, continuous function, and let $\HH$ be a compact convex set of the form~\eqref{explicit contrained set 2}. For the algorithm \eqref{eq:subspace-algo}, assume that step-sizes $\{\gamma^k\}$ satisfies \eqref{stepsize} and that the sequence $\{N^k\}$ is such that there exists $B_{\mathrm{traj}} \in \real$ satisfying $\norm{h^k} \le B_{\mathrm{traj}}$ and also $\{\widehat{\beta}^{N_k}\}$ is bounded with probability one. Then, for any $\epsilon > 0$, there exist $c_{\epsilon}(B_{\mathrm{traj}}) > 0$ and $N_{\epsilon}(B_{\mathrm{traj}}) \in \mathbb{N}$ such that $c \geq c_{\epsilon}(B_{\mathrm{traj}})$ and ${N_k \geq N_{\epsilon}(B_{\mathrm{traj}})}$ for all $k$ imply that the iterates of~\eqref{eq:subspace-algo} satisfy, with probability one, 
	\begin{equation*}
		\lim_{k \rightarrow \infty} \norm{h^k - h^*} \leq \epsilon.
	\end{equation*}
\end{proposition}
\begin{proof}
	First we show that $h^k \in \HHaff$ for all $k$. To see this, recall that $AL v= 0$ for any $v \in \real^n$. Using this in~\eqref{eq:subspace-algo} implies $Ah^{k+1} = Ah^k$ for all $k$. Consequently, for all $k$, we have $Ah^k = A h^0 = b$ and therefore $h^k \in \HHaff$. 
	
	Analogous to the proof of Proposition \ref{projected theorem}, it can be established that $\{h^k \}$ converges with probability one, in the sense of \eqref{eq:convergence to trajectory}, to a trajectory of the following dynamics
	\begin{equation} \label{eq:subspace-ode}
		\dot{\bh}(t) \! = \! - \! L \Bigl(  F\big(\bh(t)\big) \!  + \!  c\Big(\bh(t) \!  - \!  \Pi_{\HHineq}\big(\bh(t) \big) \Big) \!  - \!  e(t) \Bigr) ,
	\end{equation}
	with the initial state $\bh(0) \in \HHaff$. Here, $e(\cdot)$ is a uniformly bounded measurable map satisfying $\norm{e(t)} \leq B$ for all $t$. We will use the above fact to establish convergence of the sequence $\{h^k\}$ by analyzing the asymptotic stability of~\eqref{eq:subspace-ode}. Note that $A\dot{\bar{h}}(t) = 0$ for all $t$, and therefore a trajectory $\bar{h}(\cdot)$ of~\eqref{eq:subspace-ode}  satisfies $\bar{h}(t) \in \HHaff$ for all $t \geq 0$ as  $\bh(0) \in \HHaff$. Now consider the Lyapunov candidate
	\begin{align*}
		V(\bh ) = \frac{1}{2} \norm{\bh - \hs}^2,
	\end{align*}
	where $h^*$ is the unique solution of $\VI(\HH,F)$, that follows from strict monotonicity. As was the case for the previous result, we will first analyze the evolution of $V$ along~\eqref{eq:subspace-ode} when $e \equiv 0$. Therefore, we define the notation ${\map{X_{e\equiv0}}{\real^n}{\real^n}}$ to represent the right-hand side of~\eqref{eq:subspace-ode} with  $e \equiv 0$. The Lie derivative of $V$ along $X_{e\equiv0}$ is
	\begin{align}
		\nabla V (\bh)^\top X_{e \equiv 0}(\bh) = -(\bh &- h^*)^\top  L \Bigl(  F(\bh)  \notag  
		\\
		& +  c\big(\bh   -  \bhti \big) \Bigr). \label{eq:subspace-lie}
	\end{align}
	Since $\bar{h},h^* \in \HHaff$, we have $A(\bar{h} - h^*) = 0$ and so ${(\bar{h} - h^*) \in \SS}$. Consequently, for any vector $v \in \real^n$, we have
	\begin{align*}
	    (\bar{h} - h^*)^\top v &= (\bar{h} - h^*)^\top \big( \Pi_{\SS}(v) + \Pi_{\SS^\perp}(v) \big) \notag 
	    \\
	    & = (\bar{h} - h^*)^\top \Pi_{\SS}(v) = (\bar{h} - h^*)^\top L v.
	\end{align*}
	Using the above equality in~\eqref{eq:subspace-lie} gives
	\begin{align}
		\nabla V (\bh)^\top X_{e \equiv 0}(\bh) = -(\bh &- h^*)^\top \Bigl(  F(\bh)  \notag 
		\\
		& +  c\big(\bh   -  \bhti \big) \Bigr). \label{eq:Lie-d-subspace}
	\end{align}
	We first upper bound the second term on the right-hand side of the above equality. We have
	\begin{align}
	    & - c (\bh - h^*)^\top \big(\bh - \bhti \big) \notag    
	    \\
	    & = - c \Big( \bh - \bhti + \bhti - h^* \Big)^\top \big(\bh - \bhti \big) \notag
	    \\
	    & = - c \norm{ \bh - \bhti }^2 \notag 
	    \\
	    & \qquad + c \big( h^* - \bhti \big)^\top \big(\bh - \bhti \big) \le 0,  \label{eq:subspace-lie-second}
	\end{align}
	where for the inequality we have used the fact that  ${\big(\bar{h} - \bhti\big)^\top\big(h^* - \bhti \big) \leq 0}$ for any $\bh \in \real^n$ (see \cite[Thm. 3.1.1]{JBHU-CL:01}). Note that the inequality~\eqref{eq:subspace-lie-second} is strict whenever $\bar{h} \neq \bhti$. We now turn our attention towards the first term in~\eqref{eq:Lie-d-subspace}. Due to strict monotonicity of $F$ and the fact that $h^* \in \SOL(\HH,F)$, we obtain
	\begin{align}
	    - (\bar{h} - h^*)^\top F(\bar{h})  < - (\bar{h} - h^*)^\top  F(h^*) \le 0 \label{eq:subspace-strict}
	\end{align}
	whenever $\bar{h} \in \HH$ and $\bh \not = h^*$. The above inequality along with~\eqref{eq:subspace-lie-second} shows $\nabla V (\bh)^\top X_{e \equiv 0}(\bh) \le 0$ for any $\bh \in \HH$. However, recalling the approach in the proof of Proposition~\ref{projected theorem}, what we require in order to establish convergence is the existence of $\delta > 0$ such that
	\begin{align}
	    \nabla V (\bh)^\top X_{e \equiv 0}(\bh) \le -\delta \text{ for all } \bh \in \overline{\HH}_{\epsilon}, \label{eq:subspace-delta}
	\end{align}
	where $\overline{\HH}_{\epsilon} := \setdef{h \in \HHaff}{\norm{h - h^*} \geq \epsilon}$. We obtain this bound below. Note that the strict inequality~\eqref{eq:subspace-strict}  along with continuity of $F$ imply that for any $h \in \HH \setminus \{h^*\}$, there exists $\varepsilon_h > 0$ such that
	\begin{align}
	    -(\hat{h} - h^*)^\top  F(\hat{h}) < 0 \text{ for all } \hat{h} \in \BB_{\varepsilon_h}(h), \label{eq:subspace-nbd}
	\end{align}
	where we recall that $\BB_{\varepsilon_h}(h)$ is the open $\varepsilon_h$-ball centered at $h$. Now let $\HH_{\epsilon} := \HH \setminus \BB_{\epsilon}(h^*)$. Since $\HH$ is compact, so is  $\HH_{\epsilon}$. Using this property and~\eqref{eq:subspace-nbd}, we deduce that there exists $\varepsilon_0 > 0$ such that for every $h \in \HH_\epsilon$ we have
	\begin{align}
	    - (\hat{h} - h^*)^\top F(\hat{h}) < 0 \text{ for all } \hat{h} \in \BB_{\varepsilon_0}(h). \label{eq:subspace-nbd-two}
	\end{align}
	Next define 
	\begin{equation*}
	    \Delta_{\varepsilon_0} := \setdef{\bar{h}  \in  \HHaff \setminus \BB_{\epsilon}(h^*)  }{  \bh \not \in \BB_{\varepsilon_0}(\HH_\epsilon) \text{ and }  \norm{\bar{h}}  \leq  B_{\mathrm{traj}}}.
	\end{equation*}
	Here, $\BB_{\varepsilon_0}(\HH_\epsilon)$ is the open $\varepsilon_0$-ball of the set $\HH_\epsilon$ and ${B_{\mathrm{traj}} > 0}$ is used as an upper bound on any trajectory $\bh(\cdot)$ of~\eqref{eq:subspace-ode}. Note that $\Delta_{\varepsilon_0}$ is compact. Therefore,  there exists $B_F > 0$ satisfying
	\begin{equation}
	    -(\bar{h} - h^*)^\top F(\bar{h}) \leq B_F \text{ for all } \bh \in \Delta_{\varepsilon_0}. \label{eq:subspace-F-bound}
	\end{equation}
	Furthermore, by definition, if $\bh \in \Delta_{\varepsilon_0}$, then $\bh \not \in \HH$ and $\bh \in \HHaff$. Thus, $\bh \in \Delta_{\varepsilon_0}$ implies $\bh \not \in \HHineq$. That is, for such a point, the inequality~\eqref{eq:subspace-lie-second} holds strictly. This along with compactness of $\Delta_{\varepsilon_0}$ implies that there exists $B_{\Pi} > 0$ such that
	\begin{equation}\label{eq:subspace-negative-pi}
	    -(\bar{h} - h^*)^\top \big(\bh   -  \bhti \big) \leq -B_{\Pi} \text{ for all } \bar{h} \in \Delta_{\varepsilon_0}.
	\end{equation}
	Using~\eqref{eq:subspace-F-bound} and~\eqref{eq:subspace-negative-pi} in~\eqref{eq:subspace-lie-second} and setting $c > \frac{B_F}{B_{\Pi}}$ yields
	\begin{align}
	    \nabla V (\bh)^\top X_{e \equiv 0}(\bh) < 0 \text{ for all } \bh \in \Delta_{\varepsilon_0}. \label{eq:subspace-first-region}
	\end{align}
	Now consider $\bar{h}$ satisfying $\bar{h} \notin \Delta_{\epsilon_0} \cup \BB_{\epsilon}(h^*)$ and ${\norm{\bar{h}} \leq B_{\mathrm{traj}}}$. Note that such a point belongs to ${\HHaff \cap \BB_{\varepsilon_0}(\HH_\epsilon) \cap \BB_{B_{\mathrm{traj}}}(0)}$. Thus, by~\eqref{eq:subspace-nbd-two}, we have $-(\bh - h^*)^\top F(\bh) < 0$ for such a point. This fact combined with~\eqref{eq:subspace-first-region} leads us to the conclusion that
	\begin{align*}
	    \nabla V(\bh)^\top X_{e \equiv 0}(\bh) < 0 \text{ for all } \bh \in \overline{\HH}_\epsilon.
	\end{align*}
	Since the left-hand side of the above equation is a continuous function and $\overline{\HH}_\epsilon$ is compact, we deduce that~\eqref{eq:subspace-delta} holds. The rest of the proof is then analogous to the corresponding section of the proof in Proposition \ref{projected theorem}.
\end{proof}
\begin{remark}\longthmtitle{Practical considerations of~\eqref{eq:subspace-algo}}\label{equivalent to bounded problem}
	In Proposition~\ref{pr:subspace-cons}, for small values of $\epsilon$, one would require a large value of $c$ to ensure convergence. This may result in large oscillations of $h^k$ when $\gamma^k$ remains large. Such behavior can be prevented by either starting with small values of $\gamma^k$ or increasing $c$ along iterations, until it reaches a predetermined size. The result is then still valid but the convergence can only be guaranteed once $c$ reaches the required size. 
	
	We note that the required assumption of boundedness of $\{h^k\}$ can be ensured by constraining the iterates in $\{h^k\}$ to lie in a hyper-rectangle containing $\HH$ (cf. \cite[Page 40]{HJK-DSC:78}). However, on the boundary of the hyper-rectangle, one would have to make use of steps of the form \eqref{algorithm with error term} to ensure that the iterates remain in the feasible set. \oprocend
\end{remark} 
\subsection{Multiplier-driven algorithm}
Algorithms \eqref{algorithm with error term} and \eqref{eq:subspace-algo} involve projection onto $\HH$ or $\HHineq$ at each iteration, which can be computationally burdensome. Our next algorithm overcomes this limitation. We assume $\HH$ to be of the form \eqref{explicit contrained set 2} and introduce a multiplier variable ${\lambda \in \realnonnegative^s}$ that enforces satisfaction of the inequality constraint as the algorithm progresses. In order to simplify the coming equations we introduce the notation ${H(h,\lambda) := F(h) + Dq(h)^\top \lambda}$, where $Dq(h)$ is the Jacobian of $q$ at $h$. The \emph{multiplier-driven algorithm} is now given as
\begin{equation} \label{Lagrangian algorithm}
    \begin{split}
        h^{k + 1} &= h^k  - \gamma^k L \big(H(h^k,\lambda^k) + \widehat{\beta}^{N_k}\big), 
        \\
        \lambda^{k+1} &= \big[\lambda^{k} + \gamma^k q(h^k) \big]^+.
    \end{split}
\end{equation}
Here $L$ is as defined in \eqref{eq:subspace-projection-matrix}. Also recall that $\widehat{\beta}^{N_k}$ is the error due to empirical estimation of $F$. The next result establishes the convergence properties of~\eqref{Lagrangian algorithm} to a KKT point of the VI (see Section~\ref{section preliminaries} for definitions) and so to a solution of the VI.
\begin{proposition} \longthmtitle{Convergence of the multiplier-driven algorithm~\eqref{Lagrangian algorithm}} \label{lagrange theorem}
    Let $F$, as defined in \eqref{definition elements F}, be a strictly monotone, continuous function, and let $\HH$ be a compact convex set of the form \eqref{explicit contrained set 2}, where functions $q^i$, $i \in \until{s}$, are affine. Assume that the LICQ holds for ${h^* \in \SOL(\HH, F)}$, and let $(h^*,\lambda^*, \mu^*
    )$ be an associated KKT point. 
    For algorithm \eqref{Lagrangian algorithm}, assume that  step-sizes $\{\gamma^k\}$ satisfies \eqref{stepsize} and let $\{N_k\}$ be such that $\{\widehat{\beta}^{N_k}\}$, $\{h^k\}$, and $\{\lambda^k\}$
    are bounded with probability~one.
    Then, for any $\epsilon > 0$, there exists an $N_\epsilon \in \mathbb{N}$ such that if $N_k \geq N_\epsilon$ for all $k$, then, with probability one,
    \begin{equation*}
        \lim_{k \rightarrow \infty} \norm{h^k - h^*} \leq \epsilon.
    \end{equation*}
\end{proposition}
\begin{proof}
    Analogous to the proof of Proposition~\ref{projected theorem}, the first step establishes convergence with probability one of the sequence $\big\{(h^k, \lambda^k)\big\}$, in the sense of \eqref{eq:convergence to trajectory}, to a trajectory $\big(\bar{h}(\cdot),\bar{\lambda}(\cdot) \big)$ of the following dynamics
    \begin{subequations}\label{eq:ODE-h-lm}
        \begin{align}
            \dot{\bar{h}}(t) &= -L \Bigl( H\big(\bar{h}(t), \bar{\lambda}(t)
            \big)  + e(t) \Bigr), \label{ODE lambda part0}  
            \\
            \dot{\bar{\lambda}}(t) &= \Big[q\big(\bar{h}(t)\big)\Big]_{\bar{\lambda}(t)}^+, \label{ODE lambda part}
        \end{align}
    \end{subequations}
    with initial condition $\bh(0) \in \real^n$ and $\bar{\lambda}(0) \in \realnonnegative^l$. Note that due to the presence of this operator in~\eqref{ODE lambda part}, $\bar{\lambda}$ is contained in the nonnegative orthant along any trajectory of the system.
    The map $\bar{e}(\cdot)$ is uniformly bounded and so, as before, we have $\norm{e(t)} \le B_e$ for all $t$. The proof of convergence of the iterates to a continuous trajectory is similar to that of \cite[Theorem 5.2.2]{HJK-DSC:78} and is not repeated here for brevity. Note that, as was the case for Proposition \ref{pr:subspace-cons}, multiplication with the matrix $L$ ensure that $h^k,\bar{h}(t) \in \HHaff$ for all $k$ and $t \geq 0$. Next, we analyze the convergence of~\eqref{eq:ODE-h-lm}. We will occasionally use $\bar{x}$ as shorthand for $(\bar{h},\bar{\lambda})$. Define the candidate Lyapunov function
    \begin{equation} \label{lyapunov function}
        V(\bar{h},\bar{\lambda}) := \frac{1}{2}\big(\norm{\bar{h} - h^*}^2 + \norm{\bar{\lambda} - \lambda^*}^2\big),
    \end{equation}
    where $h^*$ is the unique solution of $\VI(\HH, F)$ and there exist $\mu^* \in \real^l$ such that $(h^*,\lambda^*, \mu^*)$ is an associated KKT point.We analyze the evolution of \eqref{lyapunov function} for the case $e \equiv 0$. Denoting the right-hand side of~\eqref{eq:ODE-h-lm} for this case by $X_{e \equiv 0}$, the Lie derivative of $V$ along~\eqref{eq:ODE-h-lm} is
    \begin{equation} \label{lagrange proof vdot}
        \begin{split}
            \nabla V(\bar{x})^\top   X_{e \equiv 0}&(\bar{x})  = -(\bar{h} - h^*)^\top H(\bar{h},\bar{\lambda})    
            \\ 
            \! \!&+ (\bar{\lambda}  \! - \! \lambda^*)^\top \! \big(q(\bar{h}) \! + \! [q(\bar{h})]_{\bar{\lambda}}^+ \! - \! q(\bar{h})\big) \!. 
        \end{split}
    \end{equation}
    Here we have dropped the matrix $L$ from the term ${(\bar{h} - h^*)^\top LH(\bar{h},\bar{\lambda})}$, which is justified by the same argument used for deriving \eqref{eq:Lie-d-subspace}. Note that for any $i$, ${([q(\bar{h})]_{\bar{\lambda}}^+)_i = (q(\bar{h}))_i}$ if $\bar{\lambda}_i > 0$. Also, if $\bar{\lambda}_i = 0$, then ${\bar{\lambda}_i - \lambda^*_i \le 0}$. Thus, we find $(\bar{\lambda} \! - \! \lambda^*)^\top([q(\bar{h})]_{\bar{\lambda}}^+ \! - \! q(\bar{h})) \leq 0$. Since $q$ is affine, we have $Dq(\bar{h}) = Dq(h^*)$ for all $\bar{h} \in \real^n$. Combined with strict monotonicity this gives, for $\bar{h} \not = h^*$,
    \begin{equation} \label{eq:derivation-in-multiplier}
        \begin{split}
            0&<(\bar{h} - h^*)^\top \big(H(\bar{h},\bar{\lambda})
            - H(h^*,\bar{\lambda}) \big)    
            \\
            &=(\bar{h} - h^*)^\top \big(H(\bar{h},\bar{\lambda})
             - H(h^*,\lambda^*) 
            \\ 
            &+ Dq(h^*)^\top\lambda^* - Dq(h^*)^\top\bar{\lambda}\big).
        \end{split}
    \end{equation}
    From \eqref{KKT} we have $-H(h^*,\lambda^*) = A^\top \mu^*$. Since we have $\bar{h},h^* \in \HHaff$ it follows that $-(\bar{h} - h^*)^\top H(h^*,\lambda^*) = 0$. Then, using the assumption that $q$ is affine, \eqref{eq:derivation-in-multiplier} gives us $- (\bar{h}  -  h^*)^\top H(\bar{h},\bar{\lambda}) < (\lambda^* - \bar{\lambda})^\top \big(q(\bar{h}) - q(h^*)\big)$ Combining these derivations, and writing $W(\bar{h})$ for the right-hand side of \eqref{lagrange proof vdot} we get that for ${\bar{h} \not = h^*}$, ${W(\bar{h}) <  (\bar{\lambda} -\lambda^{*} )^\top  q(h^*)}$. From \eqref{KKT} we have $\lambda^{*\top }q(h^*) = 0$ and $\bar{\lambda}^\top q(h^*) \leq 0$, which then implies $\nabla V(\bar{h},\bar{\lambda})  X_{e \equiv 0}(\bar{h},\bar{\lambda}) < 0$ for almost all $t$ with $\bar{h}(t) \not = h^*$. The rest of the proof is analogous to the corresponding section of the proof of Proposition \ref{projected theorem}.
\end{proof}
\begin{remark} \longthmtitle{Implementation aspect of Proposition~\ref{lagrange theorem}} \label{remark generalization}
    In Proposition \ref{lagrange theorem} we require boundedness of $\{h^k\}$, $\{\lambda^k\}$. When upper bounds on $\norm{\lambda^*}$ are known beforehand, projection onto hyper-rectangles can ensure boundedness of $\{\lambda^k\}$, while the result remains valid, (cf.~\cite[Page 40,~Theorem 5.2.2]{HJK-DSC:78}). For boundedness of $\{h^k\}$, see Remark \ref{equivalent to bounded problem}.
    \oprocend
 \end{remark}
\section{Simulations} \label{section simulations}
\begin{figure}     
    \centering
    \includegraphics[width = 0.9\columnwidth]{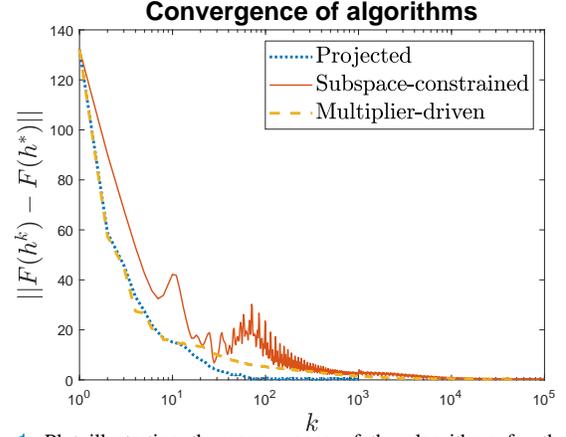}
    \vspace*{-2ex}
    \caption{\footnotesize Plot illustrating the convergence of the algorithms for the routing example explained in Section~\ref{section simulations}. The initial condition for all algorithms is set as $h_0$ defined as $(h_0)_i = 30$ for $i \in \{1,\dots,10\}$, $(h_0)_i = 60$ for $i \in \{11, \dots, 20\}$ and $(h_0)_i = 20$ for $i \in \{21,\dots, 30\}$.}
    \label{figure logscale}
\end{figure}
Here we demonstrate an application of the presented stochastic approximation algorithms for finding the solutions of a CVaR-based variational inequality. The example is an instance of a CVaR-based routing game (see Section~\ref{section routing game}) based on the Sioux Falls network~\cite{120820}. The network consists of $24$ nodes and $76$ edges. To each of the edges, we associate an affine cost function given by $C_e(f_e,u_e) = t_e(1 + u_{e} \frac{100}{c_e}f_e)$, where $f_e$ is the flow over edge $e$, and $t_e$ and $c_e$ are the free-flow travel time and capacity of edge $e$, respectively, as obtained from~\cite{120820}. The uncertainty $u_e$ has the uniform distribution over the interval $[0,0.5]$ for all edges connected to the vertices $10$, $16$, or $17$. For the rest of the edges, $u_e$ is set to zero. This defines the cost functions for all edges, and consequently defines the costs of all paths through the network as well. We consider three origin destination(OD) pairs $\WW = \{(1,19),(13,8),(12,18)\}$, and for each of these paths we select the ten paths that have the smallest free-flow travel time associated to them. The set of these 30 paths we denote as $\PP$. The demands for each OD-pair are given by $d_{1,19} = 300$, $d_{13,8} = 600$, $d_{12,18} = 200$. We aim to find a CVaR-based Wardrop equilibrium which is equivalent to finding a solution of the VI problem defined by a map ${F(h) := A h + b + \CVaR_\alpha[\xi]}$, and a feasible set $\HH = \setdef{{h \in \real^{30}}}{{h \geq 0},  {\sum_{i = 1}^{10}h_i = 300},  {\sum_{i = 11}^{20}h_i = 600}, \newline  {\sum_{i = 21}^{30}h_i = 200}}$. Here, $h,b \in \real^{30}$, ${A \in \real^{30 \times 30}}$, and ${\alpha = 0.05}$. The exact values of $A$ and $b$ and the distribution of $\xi$ are constructed using the cost functions and the network structure, see~\cite[Section 6]{AC:22} for details.

In Figure \ref{figure logscale}, we see the evolution of the error for each of the different algorithms. The stepsize sequence for the projected, subspace-constrained, and multiplier-driven algorithms are $\gamma^k = \frac{100}{100 + k}$, $\gamma^k = \frac{200}{200 + k}$, and ${\gamma^k = \min(\frac{100}{100 + k},\frac{1}{2})}$, respectively.  In addition, for the subspace-constrained algorithm we initially let $c$ depend on $k$, to prevent unstable behaviour. We used $c = \min(\frac{1}{\gamma^k},200)$. For the multiplier driven algorithm, for similar reasons, we used a modified step-size sequence for updating the multipliers $\lambda$ given by $\gamma^k_{\lambda} = 2 \gamma^k$ for $k < 1000$ and $\gamma^k_{\lambda} = 0.5 \gamma^k$, otherwise. The figure shows that all algorithms converge to a neighbourhood of the solution of the variational inequality, albeit requiring a different number of iterations. Specifically, the number of iterations taken by the projected algorithm to converge is two orders of magnitude less than that of the subspace-constrained and multiplier-driven algorithms. The quality of convergence is summarized in Table \ref{tab:mean}, where we can see both the accuracy of the achieved convergence as well as the effect of increasing the sample sizes. It is important to note that the errors shown in Fig.~\ref{figure logscale} and Table~\ref{tab:mean} are in terms of the deviation in the value of the map $\norm{F(h^k) - F(h^*)}$, rather than deviation in the solution $\norm{h^k - h^*}$. This is because the solution $h^*$ is not unique for the formulated VI.  However, for any two solutions ${h^*,\tilde{h}^* \in \SOL(F,\HH)}$ we do have $F(\tilde{h}^*) = F(h^*)$.
\begin{table}[] 
	\centering
	\begin{tabular}{l|l|l|l}
		\text{Samples per iteration}&25  &50  &100  \\ \hline
		\text{Projected}&0.3875 &0.2062  &0.1157  \\ \hline
		\text{Subspace constrained}&0.3780 &0.2015  &0.1332  \\ \hline
		\text{Multiplier-driven}&0.3889  &0.1987  & 0.1064
	\end{tabular}
    \caption{The average error of the iterates for each of the algorithms after an upper bound on the error of $0.6$, $0.3$ and $0.15$ has been achieved using $25$, $50$ and $100$ samples in each iteration respectively. The number of iterates used are $1000$, $50000$ and $100000$ for the projected, subspace-constrained and multiplier-driven algorithms respectively.}
    \label{tab:mean}
\end{table}
\section{Conclusions} \label{section conclusion}
We have considered variational inequalities defined by the CVaR of cost functions and provided stochastic approximation algorithms for solving them. We have analyzed the asymptotic convergence of these algorithms when, at each iteration, only a finite number of samples are used to estimate the CVaR. Future work will focus on analyzing the finite-time properties of the introduced algorithms and the sample complexity for a desired error tolerance of the last iterate. Finally, we wish to explore accelerated methods to solve the problem. 

\bibliography{bibfile}

\begin{thebibliography}{10}

\bibitem{FF-JSP:03}
F.~Facchinei and J.~S. Pang, {\em Finite-Dimensional Variational Inequalities
  and Complementarity Problems}.
\newblock Springer-Verlag New York, 2003.

\bibitem{UVS:13}
U.~V. Shanbhag, ``Stochastic variational inequality problems: Applications,
  analysis, and algorithms,'' {\em INFORMS TutORials in Operations Research},
  pp.~71--107, 2013.

\bibitem{YC-GL-YO:14}
Y.~Chen, G.~Lan, and Y.~Ouyang, ``Accelerated schemes for a class of
  variational inequalities,'' {\em Mathematical Programming}, vol.~165,
  p.~113–149, 2014.

\bibitem{RPR-NDS:10}
R.~Russo and N.~Shyamalkumar, ``Bounds for the bias of the empirical {CTE},''
  {\em Insurance: Mathematics and Economics}, vol.~47, pp.~352--357, 2010.

\bibitem{UR:14}
U.~Ravat, {\em On the analysis of stochastic optimization and variational
  inequality problems}.
\newblock PhD thesis, University of Illinois at Urbana-Champaign, 2014.

\bibitem{AC:19-cdc}
A.~Cherukuri, ``Sample average approximation of {CV}a{R}-based wardrop
  equilibrium in routing under uncertain costs,'' in {\em {IEEE} Conference on
  Decision and Control}, pp.~3164--3169, 2019.

\bibitem{AT-YC-MG-SM:17}
A.~Tamar, Y.~Chow, M.~Ghavamzadeh, and S.~Mannor, ``Sequential decision making
  with coherent risk,'' {\em IEEE Transactions on Automatic Control}, vol.~62,
  no.~7, pp.~332--3338, 2017.

\bibitem{CJ-LAP-FM-SM-CS:18}
C.~Jie, L.~Prashanth, M.~Fu, S.~Marcus, and C.~Szepesvári, ``Stochastic
  optimization in a cumulative prospect theory framework,'' {\em IEEE
  Transactions on Automatic Control}, vol.~42, no.~4, pp.~2867--2882, 2018.

\bibitem{HJK-DSC:78}
H.~J. Kushner and D.~S. Clark, {\em Stochastic approximation methods for
  constrained and unconstrained systems}.
\newblock Appl. Math. Sci., New York, NY: Springer, 1978.

\bibitem{VB:08}
V.~Borkar, {\em Stochastic Approximation: A Dynamical Systems Viewpoint}.
\newblock Hindustan Book Agency, 2008.

\bibitem{JV-AC:20}
J.~Verbree and A.~Cherukuri, ``Stochastic approximation for cvar-based
  variational inequalities,'' in {\em 2020 59th IEEE Conference on Decision and
  Control (CDC)}, pp.~2216--2221, IEEE, 2020.

\bibitem{DS:18}
D.~Steck, {\em Lagrange multiplier methods for contraint optimization and
  variational problems in Banach spaces}.
\newblock PhD thesis, Julius Maximilians Universität Würzburg, 2018.

\bibitem{UF-WK-GS:10}
U.~Faigle, W.~Kern, and G.~Still, {\em Algorithmic Principles of Mathematical
  Programming}, vol.~24.
\newblock Dordrecht, The Netherlands: Kluwer Academic Publishers, 2010.

\bibitem{AS-DD-AR:14}
A.~Shapiro, D.~Dentcheva, and A.~Ruszczy\'{n}ski, {\em Lectures on stochastic
  programming}.
\newblock Philadelphia, PA: SIAM, 2014.

\bibitem{MS-AC-NM:22-auto}
M.~Shakarami, A.~Cherukuri, and N.~Monshizadeh, ``Steering the aggregative
  behavior of noncooperative agents: a nudge framework,'' {\em Automatica},
  vol.~136, p.~110003, 2022.

\bibitem{JBHU-CL:93}
J.-B. Hiriart-Urruty and C.~Lemar\'echal, {\em Convex Analysis and Minimization
  Algorithms I}.
\newblock Grundlehren Text Editions, New York: Springer, 1993.

\bibitem{DCL-SRL-JJM:16}
D.~C. Lay, S.~R. Lay, and J.~J. McDonald, {\em Linear algebra and its
  applications, fifth edition}.
\newblock Pearson, 2016.

\bibitem{CDM:00}
C.~D. Meyer, {\em Matrix analysis and applied linear algebra}, vol.~71.
\newblock Siam, 2000.

\bibitem{JBHU-CL:01}
J.~Hiriart-Urruty and C.~Lemaréchal, {\em Fundamentals of Convex Analysis}.
\newblock Springer-Verlag Berlin Heidelberg New York, 2001.

\bibitem{120820}
{Transportation Networks for Research Core Team.}, ``Transportation networks
  for research.''
\newblock \url{https://github.com/bstabler/TransportationNetworks.} Accessed on
  12.08.2020.

\bibitem{AC:22}
A.~Cherukuri, ``Sample average approximation of conditional value-at-risk based
  variational inequalities,'' 2022.
\newblock Available at~\url{arxiv.org/abs/2208.11403}.

\end{thebibliography}
\bibliographystyle{ieeetr}

\end{document}